\documentclass[a4paper,11pt]{amsart}
\usepackage{amssymb}

\newtheorem{theorem}{Theorem}[section]
\newtheorem{lemma}[theorem]{Lemma}

\newtheorem{corollary}[theorem]{Corollary}

\theoremstyle{definition}

\theoremstyle{remark}
\newtheorem{remark}[theorem]{Remark}

\numberwithin{equation}{section}

\newcommand{\DD}{{\mathbb D}}

\newcommand{\NN}{{\mathbb N}}
\newcommand{\QQ}{{\mathbb Q}}

\newcommand{\CC}{{\mathbb C}}

\newcommand{\eps}{\varepsilon}

\DeclareMathOperator{\psh}{PSH}

\DeclareMathOperator{\dist}{dist}

\renewcommand{\phi}{\varphi}

\newcommand{\abs}[1]{\lvert#1\rvert}

\begin{document}

\title{Graphs that are not complete pluripolar}

\author{Armen Edigarian}

\address{Institute of Mathematics, Jagiellonian University,
Reymonta 4/526, 30-059 Krak\'ow, Poland}
\email{edigaria@im.uj.edu.pl}
\thanks{The first author was supported in part by the KBN grant
No.~5 P03A 033 21. The first author is a fellow of the 
A.~Krzy\.zanowski Foundation (Jagiellonian University)}

\author{Jan Wiegerinck}

\address{Faculty of Mathematics, University of Amsterdam,
Plantage Muidergracht 24, 1018 TV, Amsterdam, The Netherlands}
\email{janwieg@wins.uva.nl}

\date{14 February 2002}




\begin{abstract} 
Let $D_1\subset D_2$ be domains in $\CC$. Under very mild conditions
on $D_2$ we show that there exist
holomorphic functions $f$, defined on $D_1$
with the property that $f$ is nowhere extendible across
$\partial D_1$, while the graph of $f$ over $D_1$ is 
{\bf not} complete pluripolar in $D_2\times\CC$. 
This refutes a conjecture of Levenberg, Martin and Poletsky \cite{LMP}.
\end{abstract}

\maketitle


\section{Introduction}

Levenberg, Martin and Poletsky \cite{LMP} have conjectured that if $f$ is
a holomorphic function, which is defined on its maximal set of existence
$D\subset \CC$, then the graph
\begin{equation*}
\Gamma_f=\{(z,f(z)): z\in D\}
\end{equation*}
of $f$ over $D$ is a complete pluripolar subset of $\CC^2$. I.e. there exists 
a plurisubharmonic function on $\CC^2$ such that it equals $-\infty$ precisely 
on $\Gamma_f$ (see e.g.~\cite{K}).
They gave support for this conjecture in
the sense that they could prove it for some lacunary series. More support was
provided by Levenberg
and Poletsky \cite{LP} and by the second author \cite{W1,W2,W3}.
Nevertheless, in this paper we show that the conjecture is false. 

In fact we have
\begin{theorem}\label{thm1}
Let $D_1\subset D_2$ be domains in $\CC$. Assume that $D_2\setminus D_1$
has a density point in $D_2$. Then there exists a holomorphic
function $f$ with domains of existence $D_1$ such that
the graph $\Gamma_f$ of $f$ over $D_1$
is not complete pluripolar in $D_2\times\CC$.
\end{theorem}
In case $D_2\setminus D_1$ has no density point in $D_2$, it is known that
$\Gamma_f$ is complete pluripolar in $D_2\times\CC$ (see \cite{W2}).
If we take in Theorem~\ref{thm1}  for $D_1$ the unit disc $\DD$ and for
$D_2$ the whole plane $\CC$, we obtain the following corollary.
\begin{corollary}\label{cor1} There exists a holomorphic function $f$ defined
on $\DD$, which does not extend holomorphically across $\partial\DD$, such that
$\Gamma_f$ is not complete pluripolar in $\CC^2$.
\end{corollary}
Theorem \ref{thm2} then states that such a function can even be smooth 
up to the boundary of $\DD$.

The first named author thanks Marek Jarnicki, Witold Jarnicki, and
Peter Pflug for very helpful discussions.


\section{Graphs with non-trivial pluripolar hull}

The \emph{pluripolar hull} of a pluripolar set $K\subset\Omega$ is the set
\begin{equation*}
K_\Omega^\ast=\{z\in \Omega: u|_K=-\infty, u\in \psh(\Omega)
\Longrightarrow u(z)=-\infty\}.
\end{equation*}

By $\omega(z, E, D)$ we denote as usual the harmonic measure of a subset $E$ of the
boundary of a domain $D$ in $\CC$ at the point $z$ in $D$ (see e.g.~\cite{R}).

Let $A=\{a_n\}_{n=1}^\infty$ be a countable dense subset of $\partial D_1$.
Under the assumptions of Theorem~\ref{thm1}, there exists
an $a\in(\partial D_1)\cap D_2$ such that $a\in \overline{A\setminus\{a\}}$.
We may assume that $a\not\in A$.
Our function $f$ will be of the form 
\begin{equation}\label{vb}
f(z)=\sum_{j=1}^\infty \frac{c_j}{z-a_j}.
\end{equation}
We will choose $c_n$ very rapidly decreasing to $0$. 
In particular, 
\begin{equation}\label{eq1}
\sum_{j=1}^\infty|c_j|<+\infty,
\end{equation}
so the limit in \eqref{vb} exists and is a holomorphic function on $D_1$.

Moreover, we will choose $c_n$ such that
$\sum_{n=1}^\infty\frac{|c_j|}{|a-a_j|}<+\infty$. Hence,
the series (\ref{vb}) will converge at $z=a$.
We will denote its limit by $f(a)$. We will prove a
version of Theorem~\ref{thm1} that elaborates on \eqref{vb}. However,
no statement about extendibility is made at this point.
\begin{theorem}\label{example} Let $D_1\subset D_2$ be domains in $\CC$, such
that $D_2\setminus D_1$ has a density point in $D_2$.
There exists a sequence $\{R_n\}_{n=1}^\infty$ of positive numbers such that 
for any sequence of complex numbers $\{c_n\}_{n=1}^\infty$ with
$|c_n|\le R_n$ we have $(a,f(a))\in (\Gamma_f)^\ast_{D_2\times\CC}$,
where $f$ is given by \eqref{vb}.
Here $\Gamma_f$ is the graph of $f$ over $D_1$.
\end{theorem}

\begin{proof} 
We may assume that $a=0$. For $b\in\CC$ and $r>0$ we set
$\DD(b,r)=\{z\in\CC: \abs{z-b}<r\}$, $\DD_r=\DD(0,r)$
and $\DD=\DD_1$, the unit disc.
Put $$\pi_{n}^k(z)=e^{\frac{2\pi ki}{n}}z,\quad 
z\in\CC,\quad k,n\in\NN,$$ and $$B=\cup_{n=1}^\infty\cup_{k=1}^n\pi_{n}^k(A).$$ 
Note that $0\not\in B$ and that $B$ is also countable (and therefore thin at 
$0$). By Corollary 4.8.3 in \cite{K}, there exists an open set $U\supset B$ 
such that $U$ is thin at $0$. 

\smallskip
\emph{Step 1.} We construct a sequence of radii $\{\rho_n\}_{n=1}^\infty$
with special properties, the main one being that
$\cup_{n}\cup_{k=1}^n \pi_{n}^k (\DD(a_n,\rho_n))$ is thin at 0. 

It is a corollary of Wiener's criterion (see \cite{R}, Theorem~5.4.2)
that there exists a sequence $r_n\to0$ such that 
\begin{equation}\label{eq3}
\partial\DD_{r_n}\cap U=\varnothing,\quad n\in\NN.
\end{equation}
Since $U$ is thin at $0$, there exists a subharmonic function $u$ on $\CC$
such that 
\begin{equation*}
\limsup_{U\ni z\to 0} u(z)=-\infty< u(0)
\end{equation*}
(see e.g.~Proposition 4.8.2 in \cite{K}).
Moreover, by scaling and adding a constant, we can assume that $u(0)=-\frac12$
and $u<0$ on $\DD$. By \eqref{eq3} there exists a $\rho>0$ such that
$\overline\DD_\rho\subset D_2$, $\partial\DD_\rho\cap D_1\not=\varnothing$, 
$\partial\DD_\rho\cap U=\varnothing$, and
$u\le-1$ on $U\cap \DD_\rho$ (take $\rho=r_n$ with sufficiently big $n$).

Let $J\subset\partial\DD_\rho\cap D_1$ be a closed arc. We can assume that 
\begin{equation*}
J=\big\{e^{i\theta}\rho:\frac{2\pi k_0 }{n_0}\le\theta
\le\frac{2\pi (k_0+1) }{n_0}\big\}
\end{equation*}
for some $k_0,n_0\in\NN$.

Now we choose a sequence of positive numbers $\rho_n\in(0,1)$, $n\in\NN$,
in the following way:
\begin{enumerate}
\item Let $0<\rho_1<1$ be such that
	\begin{enumerate}
	\item $\cup_{k=1}^{n_0}\pi_{n_0}^k\big(\DD(a_1,\rho_1)\big)\subset U$;
	\item $\DD_\rho\setminus\cup_{k=1}^{n_0}
	\pi_{n_0}^k\big(\overline\DD(a_1,\frac{\rho_1}{2})\big)$ is connected.
	\end{enumerate}
\item Assume that $\rho_1,\dots,\rho_{n-1}$ are chosen. 
Choose $0<\rho_n<1$ such that 
	\begin{enumerate}
	\item $\cup_{k=1}^{n_0}\pi_{n_0}^k\big(\DD(a_n,\rho_n)\big)\subset U$;
	\item $\DD_\rho\setminus\cup_{j=1}^n\cup_{k=1}^{n_0}
	\pi_{n_0}^k\big(\overline\DD(a_j,\frac{\rho_j}{2})\big)$ is connected.
	\end{enumerate}
\end{enumerate}
Put $Y_n=\bigcup_{j=1}^n\bigcup_{k=1}^{n_0}
\pi_{n_0}^k\big(\overline\DD(a_j,\frac{\rho_j}{2})\big)$.
So, $Y_n\subset U$ is a closed set such that $\DD_\rho\setminus Y_n$
is a domain and $\partial\DD_\rho\cap Y_n=\varnothing$ for any $n\in\NN$.

\smallskip
\emph{Step 2.} We want to show that
\begin{equation}\label{eq4}
\omega(0,\partial\DD_\rho,\DD_\rho\setminus Y_n)\ge\frac12,\quad n\in\NN.
\end{equation}
Fix $n\in\NN$.
Put $v_n(z)=-\omega(z,\partial\DD_\rho,\DD_\rho\setminus Y_n)+u(z)$.
It suffices to show that 
\begin{equation}\label{eq6}
v_n\le-1\quad\text{ on }\DD_\rho\setminus Y_n.
\end{equation}
Observe that 
$-\omega(\cdot,\partial\DD_\rho,\DD_\rho\setminus Y_n)\le0$ and $u(\cdot)\le0$
on $\DD_\rho\setminus Y_n$.
Moreover, we have $\limsup_{z\to\partial\DD_\rho} 
-\omega(z,\partial\DD_\rho,\DD_\rho\setminus Y_n)\le-1$ and
$\limsup_{z\to Y_n} u(z)\le -1$. So, from the maximum principle for
the subharmonic function $v_n$ we get \eqref{eq6} and, therefore, \eqref{eq4}.

\smallskip
\emph{Step 3.} Here we want show that
\begin{equation}\label{eq5}
\omega\big(0,J,\DD_\rho\setminus
\cup_{j=1}^n\overline\DD(a_j,\frac{\rho_j}2)\big)\ge\frac1{2n_0},\quad n\in\NN.
\end{equation}
Put
\begin{equation*}
w_n(z)=\omega\big(z,\partial\DD_\rho,\DD_\rho\setminus Y_n\big)-
\sum_{k=1}^{n_0}\omega\big(z,\pi_{n_0}^k (J),\DD_\rho\setminus Y_n\big).
\end{equation*}
Note that $\cup_{k=1}^{n_0} \pi_{n_0}^k(J)=\partial\DD_\rho$.
Again from the maximum principle we obtain that $w_n\le0$ on 
$\DD_\rho\setminus Y_n$, $n\in\NN$.

Because $\pi_{n_0}^k(\DD_\rho\setminus Y_n)=\DD_\rho\setminus Y_n$,
for any $k,n\in\NN$, we find 
\begin{equation*}
\omega\big(0,\pi_{n_0}^k (J),\DD_\rho\setminus Y_n\big)=
\omega\big(0,J,\DD_\rho\setminus Y_n\big),\quad k\in\NN.
\end{equation*}

Hence, 
\begin{equation*}
\omega\big(0,J,\DD_\rho\setminus
\cup_{j=1}^n\overline\DD(a_j,\frac{\rho_j}2)\big)\ge
\omega\big(0,J,\DD_\rho\setminus Y_n\big)\ge\frac{1}{2n_0},\quad n\in\NN.
\end{equation*}

\smallskip
\emph{Step 4.} 
Let $\{R_n\}_{n=1}^\infty$ be a sequence of positive numbers
such that $C_1:=\sum_{n=1}^\infty\frac{R_n}{\rho_n}<+\infty$ and, therefore,
$\sum_{n=1}^\infty R_n<C_1$
(take e.g.~$R_n=\frac{\rho_n}{n^2}$). Consider any sequence
of complex numbers $\{c_n\}_{n=1}^\infty$ with $|c_n|\le R_n$ and 
let $f$ be defined by \eqref{vb}.

Put 
$$f_n(z)=\sum_{j=1}^{n}\frac{c_j}{z-a_j}-\sum_{j=n+1}^\infty\frac{c_j}{a_j},
\quad n\in\NN.$$
Then $|f_n(z)|\le 2C_1$ for every $z\in\DD_\rho\setminus
\cup_{j=1}^n\overline{\DD}(a_j,\frac{\rho_j}{2})$ and all $n$.

Let  $h\in\psh(D_2\times\CC)$ have the property that 
$h\big(z,f(z)\big)=-\infty$, $z\in D_1$. The function $s_n$ defined
on $D_2\setminus\{a_1,\dots,a_n\}$ by
$s_n(z):=h(z,f_n(z))$ is subharmonic.
Let $A_n:=\sup_{z\in J} s_n(z)$ and let
$C_2:=\sup_{z\in\overline\DD_\rho,|w|\le 2C_1} h(z,w)$.
Then $A_n\to\sup_{z\in J} h(z,f(z))=-\infty$ as $n\to\infty$.

From the two-constant theorem (see e.g.~\cite{R}, Theorem~4.3.7)
we infer
\begin{equation*}
\frac{C_2-s_n(0)}{C_2-A_n}\ge
\omega\big(0,J,
\DD_\rho\setminus\cup_{j=1}^n \overline\DD(a_j,\frac{\rho_j}{2})\big)
\ge\frac{1}{2n_0},\quad n\in\NN.
\end{equation*}
Letting $n\to\infty$, we conclude that
 $h(0,f(0))=s_n(0)=-\infty$ and therefore
$(0,f(0))\in(\Gamma_f)^\ast_{D_2\times\CC}$.
\end{proof}

For the proof of Theorem~\ref{thm1} we need to know that the function
defined by \eqref{vb} is not extendible across the boundary of $D_1$.
This will be done in the next section.

\section{Non-extendible sums}

Without additional conditions on $a_n$ and $c_n$ a function defined
by \eqref{vb} may well  extend holomorphically beyond the boundary of $D_1$,
think of Lambert-type series
$\sum_{n=1}^\infty c_n\frac{z^n}{1-z^n}$, cf.~\cite{Knopp}.
It may  even yield 0 on $D_1$,
cf.~\cite{BSZ}. We will see that suitable choice of $a_n$ and $c_n$ prevents
 this
from happening. We are grateful to Marek~Jarnicki and
Peter~Pflug who suggested the idea of the proof of the next lemma.

\begin{lemma}\label{lem1}
Let $D$ be a domain in $\CC$. Then there exist a dense subset 
$A=\{a_n\}_{n=1}^\infty$ of $\partial D$ and a sequence 
$\{R_n\}_{n=1}^\infty$ of positive numbers such that for any
sequence of complex numbers $\{c_n\}_{n=1}^\infty$ with
$0<|c_n|\le R_n$ the holomorphic function $f$ given by \eqref{vb}
is not holomorphically extendible across $\partial D$.
\end{lemma}

\begin{proof} 
Let $B=\{b_n\}_{n=1}^\infty$ be a dense subset of $D$ 
(take e.g. $U=D\cap\QQ^2$).
For any $b_n\in B$ there exists a point $a\in\partial D$
such that $\dist(b_n,\partial D)=|b_n-a|$. We denote by $a_n$ one of them.
Set $A=\{a_n\}_{n=1}^\infty$. Note that $A$ is a dense subset of $\partial D$.
Taking subsequence of $\{a_n\}_{n=1}^\infty$ we may assume that
$a_i\ne a_j$, $i\ne j$.

Fix $n\in\NN$. Let $B_{n}=\{z\in D:\dist(z,\partial D)=|z-a_n|\}\subset D$.
Note that $B_n\cup\{a_n\}$ is a closed set on the plane and 
$\widetilde B=\cup_{n=1}^\infty B_n$ is dense in $D$ (because 
$\widetilde B\supset B$). Moreover, if $z_0\in B_n$ then
the open segment with the ends
at the points $z_0$ and $a_n$ is contained in $B_n$. 

For any $j\in\NN$ we put $\epsilon_{nj}=\dist(a_j,B_{n})$. 
Since $a_j\not\in B_n$ for $j\ne n$, 
we see that $\epsilon_{nj}>0$ for $j\ne n$.

Put 
\begin{equation*}
R_j=\frac{\min\{\epsilon_{1j},\dots,\epsilon_{(j-1)\, j}\}}{j^2},
\quad j\in\NN.
\end{equation*}
For any $n\in\NN$ and any $j>n$ we have $R_j\le\frac{\epsilon_{nj}}{j^2}$
and therefore 
\begin{equation*}
\sum_{j\not=n}\frac{R_j}{|z-a_j|}\le
\sum_{j\ne n}\frac{R_j}{\epsilon_{nj}}<+\infty,\quad z\in B_n.
\end{equation*}

Take a sequence of complex numbers 
$\{c_n\}_{n=1}^\infty$ with  $0<|c_n|\le R_n$.
Then for a fixed $n\in\NN$ we have
\begin{multline}\label{eq7}
\liminf_{B_n\ni z\to a_n}|(z-a_n)f(z)|\\
\ge |c_n|-\lim_{B_n\ni z\to a_n} |z-a_n|\cdot
\limsup_{B_n\ni z\to a_n}\sum_{j\not=n}
\frac{|c_j|}{|z-a_j|}=|c_n|>0.
\end{multline}

Observe that for any $n\in\NN$ the Taylor series at any point of $z_0$
$B_n$ has a radius of convergent equal to $\dist(z_0,\partial D)$
(because of \eqref{eq7} and $|z_0-a_n|=\dist(z_0,\partial D)$).
Hence, by Lemma 1.7.5 from \cite{JP} we see that $D$ is the domain
of existence of $f$.
\end{proof}

\begin{proof}[Proof of Theorem~\ref{thm1}]
If a set $E\subset\Omega$ is complete pluripolar in a domain 
$\Omega$, then $E^\ast_{\Omega}=E$. 
By Lemma~\ref{lem1} and Theorem~\ref{example} 
there exists a  holomorphic  function $f$ on $D_1$ for 
which $D_1$ is a domain of existence and 
$(\Gamma_f)_{D_2\times\CC}^\ast\not=\Gamma_f$. Hence, $\Gamma_f$ is not
complete pluripolar in $D_2\times\CC$.
\end{proof}

\begin{theorem}\label{thm2}
There exists a sequence $\{a_n\}_{n=1}^\infty\subset\CC\setminus{\overline\DD}$
and a sequence $\{c_n\}_{n=1}^\infty$ such that the function $f$
defined by \eqref{vb} is  $C^\infty$ on $\overline\DD$, is nowhere extendible
over the boundary of $\DD$, while $\Gamma_f$ is not complete pluripolar in 
$\CC^2$.
\end{theorem}

\begin{proof}
Let $r_j=1+1/(j+1)$. The sequence $a_n$ is formed by 
$$a_{2^j+k}=r_je^{2\pi i \frac{k}{2^j}},\quad k=0,\ldots 2^j-1,\ j=0,1,\ldots.$$
The proof of Theorem~\ref{example} provides us with a sequence $\{R_n\}$
such that
for every sequence $\{c_n\}$ with $|c_n|<R_n$ the series \eqref{vb} represents
a function on $\DD$,
the graph of which is not complete pluripolar. Assembling all 
$a_n\in C(0,r_j)$ we find that there exists a sequence $\{R'_j\}$ such that for
every choice of $0<\eps_j<R_j$ the function $f_\eps$ on $\DD$ defined by
\begin{equation}
f_\eps(z)=\sum_{j=0}^\infty\frac{\eps_j}{r_j^{2^j}-z^{2^j}}
\label{series}\end{equation}
has a graph that is not complete pluripolar.

We observe that independently of the choice of $\eps_j$
$$\sum_{j=1}^{n}\frac{\eps_j}{r_j^{2^j}-z^{2^j}}
=\sum_{k=1}^\infty d_{n,k} z^k$$
is holomorphic on $\DD_{r_n}$ with singularities on the boundary. 
(The $d_{n,k}$ are defined by the equality.)
Therefore we have
$\limsup_{k\to\infty}|d_{n,k}^{1/k}|=1/r_{n}$.
Hence there is a $k_n$ such that
\begin{equation}|d_{n,k_n}|> r_{n+1}^{-k_n}.\label{rad1}
\end{equation}

We will now make an appropriate choice for the $\eps_j$ to insure
that $f$ cannot be extended over the boundary of $\DD$. Along the way
we will determine constants
$C_j$ that are needed for smoothness at the boundary.

Choose $\eps_0=R'_0$.
Then $$f_{\eps_0}(z)=\frac1{r_0-z}=\sum_{k=1}^\infty d_{0,k}z^k$$ with
$\limsup_{k\to\infty}|d_{0,k}^{1/k}|=1/r_{0}$; in particular there is a $C_0$
such that 
$|d_{0,k}|<{C_0}$.

Suppose $\eps_0,\ldots,\eps_{n-1}$ and $C_0,\ldots,C_{n-1}$
have been chosen in such a way that we have found $k_0,\ldots,k_{n-1}$ with
\begin{equation}
|d_{l,k_j}|>r_{j+1}^{-k_j}\quad
\text{for\ } j=0,\ldots,n-1,\ l=j,\ldots n-1\label{radius}
\end{equation}
and 
\begin{equation}
|d_{j,k}|<\frac{C_l}{k^l}\quad \text{for\ } l=0,\ldots, n-1, \ j= 0,\ldots, n-1
\quad\text{and all \ } k.
\label{smooth}
\end{equation}
Then choose  
\begin{equation}
C_n> \sup_k |d_{j,k}|k^n\quad j=0,\ldots,n-1.\label{Cn}
\end{equation} 
This is finite because of \eqref{radius}. Next choose $\eps_n<R'_n$
so small that
\begin{enumerate}
\item The inequality
\eqref{smooth} holds for $l=0,\ldots,n$ and $j= 0,\ldots, n$.
This is possible because of \eqref{radius}.  
\item 
\begin{equation*}
|d_{n,k_j}|>r_{j+1}^{-k_j}\quad
\text{for\ } j=0,\ldots,n-1,
\end{equation*}
which is again possible because of \eqref{radius}.
\end{enumerate}
Having chosen $\eps_n$,
 we can by \eqref{rad1} choose $k_n$ so large that
$|d_{n,k_n}|>r_{n+1}^{-k_n}$.

Observe that the coefficients $d_{n,k}$ converge to the coefficients
$d_k$ of the power series expansion of $f_\eps$ as $n\to \infty$.
From \eqref{radius} we see that $|d_{k_j}|> r_{j+1}^{-k_j}$ so that
the radius of convergence of the power series
of $f_\eps$ is at most 1, and since $f_\eps$ is holomorphic on $\DD$,
it equals 1. So $f_\eps$ has a singular point $b$ on $C(0,1)$.
We split $f_\eps$ as
$$f_\eps=f_1+f_2=\left(\sum_{j=0}^{n-1}+\sum_{j=n}^\infty\right)
\frac{\eps_j}{r_j^{2^j}-z^{2^j}}.$$
Then $f_1$ is holomorphic in a neighborhood of the closed unit disc and $f_2$ has at least one singular point on $C(0,1)$, but $f_2$ is invariant under rotation
over a $2^n$-th root of unity, which implies that there is a singularity in
each arc of length $>2\pi/2^n$. Therefore, $f$ can nowhere be extended 
analytically over $C(0,1)$.

Next we show that $f$ is smooth up to the boundary of $\DD$.
We have to show that there exist constants $C_l>0$ such that for every $l$
$$|d_k|\le \frac{C_l}{k^l},\quad \text{for all\ } k, $$
but this follows from \eqref{smooth}.
\end{proof}
\smallskip
\begin{remark} Let $D$ be a  domain in $\CC$ and let $A$ be a
 closed polar subset of $D$.
Using methods presented in this paper and in \cite{W2},
the authors give in \cite{AW} a complete characterization of the
holomorphic functions $f$
on $D\setminus A$ such that 
$\Gamma_f=\{(z,f(z)):z\in D\setminus A\}$ is complete pluripolar 
in $D\times\CC$.
\end{remark}
%

\bibliographystyle{amsplain}


\end{document}